\documentclass[10pt,reqno]{amsart}

\usepackage[utf8]{inputenc}
\usepackage{comment}
\usepackage{amsmath}
\usepackage{mathtools}
\usepackage{graphicx}
\usepackage{amsfonts}
\usepackage{amsthm}
\usepackage{amssymb}
\usepackage{mathrsfs}	
\numberwithin{equation}{section}
\usepackage[all]{xy}
\usepackage{color}
\usepackage[left=3cm,right=3cm]{geometry}
\usepackage{bm} %bold mathcal
\usepackage{wasysym}
\usepackage{ytableau}
\usepackage{hyperref}
\usepackage{cleveref}

\usepackage{todonotes}

\DeclareMathOperator{\ch}{ch}

\DeclareMathOperator{\Hook}{Hook}
\DeclareMathOperator{\hook}{h}

\DeclareMathOperator{\End}{End}
\DeclareMathOperator{\Ss}{\mathfrak{S}}
\newcommand{\Cl}{\mathrm{Cl}}

\newcommand{\D}{\mathsf{D}}
\newcommand{\DD}{\mathbb{D}}
\newcommand{\blambda}{\mathsf{m}}
\newcommand{\power}{\mathsf{\mathbf{p}}}
\newcommand{\Rep}{{\Lambda}}

\newcommand{\Schur}{\mathsf{\mathbf{s}}}

\newcommand{\Q}{\mathbb{Q}}
\newcommand{\Z}{\mathbb{Z}}
\newcommand{\N}{\mathbb{N}}
\newcommand{\llangle}{\langle \!\langle}
\newcommand{\rrangle}{\rangle \!\rangle}
\newtheorem{lemma}{Lemma}[section]
\newtheorem{proposition}[lemma]{Proposition}
\newtheorem{thm}[lemma]{Theorem}
\newtheorem*{thmA}{Theorem A}
\newtheorem*{thmB}{Theorem B}

\newtheorem{corollary}[lemma]{Corollary}

\theoremstyle{definition}
\newtheorem{definition}[lemma]{Definition}
\newtheorem{example}[lemma]{Example}
\newtheorem{remark}[lemma]{Remark}
\title{A plethystic chain rule}

\author{Alessandro D'Andrea}
\address{\newline Alma Mater studiorum Universit\`a di Bologna\hfill\newline Dipartimento di Matematica\hfill\newline Piazza di porta San Donato 5, 40126 Bologna, Italy}
\email[A.~D'Andrea]{a.dandrea@unibo.it}

\author{Enrico Fatighenti}
\address{\newline Alma Mater studiorum Universit\`a di Bologna\hfill\newline Dipartimento di Matematica\hfill\newline Piazza di porta San Donato 5, 40126 Bologna, Italy}
\email[E.~Fatighenti]{enrico.fatighenti@unibo.it}

\author{Claudio Onorati}
\address{\newline Università della Calabria \hfill \newline Department of Mathematics and Computer Science\hfill\newline via Pietro Bucci, 87036 Arcavacata di Rende (CS), Italy} \email{claudio.onorati@unical.it}

\begin{document}
	
	\begin{abstract}
		We consider a derivation $\D$ on the ring $\Lambda$ of symmetric functions and investigate its combinatorial, algebraic and geometric properties.
		More precisely, we show that $\D$ restricts to a quasi-isometry, with respect to the Hall product, on the graded component of $\Lambda$ of each positive degree and provide a chain-rule formula with respect to the plethysm operation.
		Furthermore, we relate the geometry of the  Schur functions supporting $\D(f)$, where $f\in \Lambda$ is a homogeneous element, to that of $f$.
	\end{abstract}
	
	
	
	\maketitle
	
	\tableofcontents
	
	%%%%%%%%%%%%%%%%%%%%%%%%%%%%%%%%%%%%%
	%%%%%%%%%%%%%%%%%%%%%%%%%%%%%%%%%%%%%
	\section*{Introduction}
	
	The ring of symmetric functions $\Rep$ is endowed with a fundamental additional operation, known as \emph{plethysm} (\cite{Littlewood}) and denoted by $f[g]$ for $f,g\in\Rep$. 
	If $\Schur_\lambda$ and $\Schur_\mu$ are Schur symmetric functions, then 
	\[ \Schur_\lambda[\Schur_\mu]=\sum_{\nu} c^\nu_{\lambda,\mu}\,\Schur_\nu, \]
	for some non-negative integer coefficients $c^\nu_{\lambda,\mu}$.
	The problem of computing the coefficients $c^\nu_{\lambda, \mu}$ has a long history, and while several solutions are known, they are often computationally inefficient. In fact, even deciding whether a given coefficient $c^\nu_{\lambda, \mu}$ is non-zero can be challenging. This difficulty has been formalized by results in \cite{FI} concerning the computational complexity of the problem. (We are grateful to an anonymous referee for bringing this connection to our attention.) In more detail, the problem of deciding whether the coefficient $c^{\nu}_{(n),(m)}$ is non-zero (for a fixed $m\geq3$) was shown to be $\mathsf{NP}$-hard (\cite[Theorem~3.5]{FI}), and it remains open whether this problem belongs to the class $\mathsf{NP}$, with experts' opinions divided. Similarly, determining the coefficients $c^{\nu}_{(n),(3)}$ is $\mathsf{\# P}$-hard (\cite[Theorem~3.6]{FI}). Negative results of a different combinatorial flavor can also be found in \cite{KM18}.
	
	The standard reference for problems in this area is \cite[Section~3]{Stanley}, and we also refer to \cite{LoehrRemmel2011,dBPW21,MysteryOfPlethysms} for more recent surveys of results. For formulas regarding special plethysm coefficients, see \cite[I.8, Example~9]{Macdonald} (as well as \cite{Boffi}).
	
	In this paper, we present an alternative recursive algorithm to compute plethysms. Our method exploits the algebraic properties of a derivation 
	\[ \D\colon\Rep\longrightarrow\Rep, \]
	defined in terms of the power sum symmetric functions $\power_n$ as
	\[ \D=\sum_{n\geq1}n\frac{\partial}{\partial\power_n}. \]
	We provide a combinatorial interpretation of the action of $\D$ on Schur functions in Proposition \ref{prop:geometric interpretation of D}. We note that while the derivation $\D$ has appeared previously in the literature (see Remark \ref{rmk:D gia nota}), it has not, to the best of our knowledge, been applied to the computation of plethysms.
	
	Our first main result establishes a relationship between the derivation $\D$ and the Hall inner product.
	
	\begin{thmA}[Theorem \ref{thm:quasi isometria}]
		For every $d>0$, the restriction of the derivation $\D$ to $\Rep^d$ is a quasi-isometry. More precisely, for any $f,g\in\Rep^d$, we have 
		\[ \langle\D(f),\D(g)\rangle=d\langle f,g\rangle. \]
	\end{thmA}
	This isometry property, combined with the fact that $\D$ is a derivation, yields a recursive procedure to compute any Littlewood-Richardson coefficient in terms of coefficients of lower degree (see Corollary~\ref{cor:LR}).
	\medskip
	
	Next, we establish an explicit formula for the action of $\D$ on a plethysm. To this end, we introduce a slight variant of the operator: let $\Rep\Psi \subset \End_\Q \Rep$ be the subring generated by the \emph{Adams operations} $\Psi^n\colon f\mapsto \power_n[f]$ and all multiplication maps by elements in $\Rep$. We define the operator
	\[ \DD\colon\Rep\longrightarrow\Rep\Psi,\qquad \DD=\sum_{n\geq1} n\frac{\partial}{\partial\power_n} \circ \Psi^n, \]
	which is essentially a reformulation of $\D$, when computed on homogeneous elements from $\Rep$.
	
	Our second main result is a chain rule formula for the plethystic operation.
	
	\begin{thmB}[Theorem \ref{thm:plethysm product}]
		For any $f,g\in\Rep$, we have
		\begin{equation*}
			\DD(f[g])=(\DD(f)[g])\circ \DD(g). 
		\end{equation*}
	\end{thmB}
	The right-hand side of this equality involves only plethysms between symmetric functions of lower degree, yielding an effective recursive method for computing $\DD(f[g])$, and consequently, $\D(f[g])$. 
	
	Since Schur functions form an orthonormal basis with respect to the Hall inner product, Theorem~A and Theorem~B can be combined to recursively compute individual plethysm coefficients $c^\nu_{\lambda, \mu}$ (see the end of Section~\ref{section:plethysms} for some examples). We emphasize that our recursion is distinct from other recursive formulas for plethysm found in the literature (see, for example, \cite{EPW,LO23,MysteryOfPlethysms}).
	\medskip
	
	For any symmetric function $f \in \Rep$, let us write $f=\sum_\lambda f_\lambda\,\Schur_\lambda$ and define the support of $f$ as the region obtained by taking the union of all Young diagrams $\lambda$ (viewed as subsets of the positive quadrant) for which $f_\lambda \neq 0$.
	In Section \ref{section:applications}, we prove Theorem \ref{thm:supported on columns}, answering the following question: to what extent does the support of $f$ coincide with that of $\D(f)$?
	
	The naive expectation that $f$ and $\D(f)$ should have the same support fails in general. Indeed, $\D(\power_n) = n = n \Schur_{\emptyset}$ is a multiple of the Schur function corresponding to the empty partition, while $\power_n$ is supported on hook partitions, which occupy a non-trivial region. However, this is essentially the only exception; Theorem \ref{thm:supported on columns} provides a precise answer for the unbounded regions defined by the first $k$ columns. Analogous statements hold for the first $k$ rows and, with minor modifications, for all rectangular regions. We expect a generalized statement to hold for regions of arbitrary shape.
	
	Theorem \ref{thm:supported on columns} complements the recursive nature of Theorem~B beautifully. We apply it to give a purely recursive proof of the well-known fact that if the partition $\mu$ is supported on $k$ columns, then every partition $\nu$ appearing in the plethysm $\Schur_\lambda[\Schur_\mu]$ is supported on $k|\lambda|$ columns (see Corollary \ref{cor:alcuni coefficienti sono zero}).
	
	While this statement is standard and follows easily from representation-theoretic considerations (e.g.\ by viewing $\Schur_\lambda[\Schur_\mu]$ as a summand of $\Schur_\mu^{\otimes |\lambda|}$), our focus is on the proof technique itself. Notably, our proof avoids any use of the Littlewood--Richardson rule. Broader generalizations of this claim can be found in \cite{PW16} (see also \cite[Theorem~1.5]{dBPW21}, \cite{KM16}). 
	
	None of the authors is a specialist algebraic combinatorialist, and we hope that these techniques find deeper applications in more expert hands. 
	The definition of the derivation $\DD$ was originally inspired by geometric problems concerning Schur functors \cite{Casimiri}, and in particular by the problem of computing the Chern character of Schur bundles, which we completely solve in \cite{DFO:Chern}.
	
	%%%%%%%%%%%%%%%%%%%%%%%%%%%%%%%%%%%%%%
	\subsection*{Plan of the paper} 
	In Section \ref{section:preliminaries} and Section \ref{section:derivation} we set our notation and recall some results from the theory of symmetric functions. In Section \ref{section:remarkable derivation}, we introduce the derivation $\D$ and prove Theorem~A; for readers more familiar with the combinatorics of representations of the symmetric group, we translate in Section \ref{section:symmetric group II} our results in this setting. In this section, we also give a recursive formula to compute the product of two symmetric functions. Theorem~B is then proved in Section \ref{section:plethysms}, where we also collect some explicit examples. In Section \ref{section:applications}, we prove Theorem \ref{thm:supported on columns} concerning the support properties of $\D$, together with some further applications. 
	
	%%%%%%%%%%%%%%%%%%%%%%%%%%%%%%%%%%%%
	\subsection*{Acknowledgments}
	We thank Michele D'Adderio for advice, comments and references. We also wish to thank Laurent Manivel for comments and criticism. We are deeply indebted to an anonymous referee, whose precious and thoughtful comments have considerably improved both the content and the presentation of this paper.
	
	This research has been partially funded by the European Union - NextGenerationEU under the
	National Recovery and Resilience Plan (PNRR) - Mission 4 Education and research - Component 2
	From research to business - Investment 1.1 Notice Prin 2022 - DD N. 104 del 2/2/2022, from title
	“Symplectic varieties: their interplay with Fano manifolds and derived categories”, proposal code
	2022PEKYBJ – CUP J53D23003840006.
	The second and third authors are members of the INDAM-GNSAGA group.

	%%%%%%%%%%%%%%%%%%%%%%%%%%%%%%%%%%%%%
	%%%%%%%%%%%%%%%%%%%%%%%%%%%%%%%%%%%%%
	\section{Preliminaries}\label{section:preliminaries}
	Here we set up the notation and recall some basic facts about symmetric functions. For a general introduction, see \cite{Zabrocki}.
	
	%%%%%%%%%%%%%%%%%%%%%%%%%%%%%%%%%%%%
	\subsection{Partitions and Young diagrams}\label{section:partitions and Young diagrams}
	
	Let $\lambda$ be an integral partition. We use the two equivalent notations
	\[ \lambda=(\lambda_1,\dots,\lambda_r)=(1^{\blambda_1}\cdots t^{\blambda_{t}}) \]
	where $\lambda_1\geq\cdots\geq\lambda_r>0$ and $\blambda_k=\sharp\{j\mid \lambda_j=k \}$ are the multiplicities. The index $r$ is the \emph{length} of the partition, also denoted by $\ell(\lambda)$. Notice the equality $r=\ell(\lambda)=\sum_j\blambda_j$. 
	
	The \emph{size} of $\lambda$ is the integer
	\[ |\lambda|=\sum_j\lambda_j=\sum_k k\,\blambda_k\,. \]
	If $\lambda$ is a partition, then its conjugate $\lambda^*$ is the partition such that $\lambda^*_k=\sharp\{j\mid \lambda_j\geq k\}$. Notice that $\blambda_k=\lambda^*_k-\lambda^*_{k+1}$.
	\medskip
	
	Given a partition $\lambda$ as above, we abuse notation and keep denoting by $\lambda$ the associated \emph{Young diagram}. Recall that this is the diagram composed by $|\lambda|$ boxes, divided in $\ell(\lambda)$ rows of lengths $\lambda_1,\dots,\lambda_{\ell(\lambda)}$. For example,
	
	\[  
	\begin{ytableau}
		~ & ~ & ~ & ~ \\
		~ & ~ & ~\\
		~  
	\end{ytableau} = (4,3,1) = (1^{1}\,2^0\,3^1\,4^1).
	\]
	When writing a partition using the notation with multiplicities, we will generally avoid writing the multiplicity $0$, thus the example above will be simply written as $(1^1\,3^1\,4^1)$.
	\medskip
	
	Taking the conjugate of a partition corresponds to flipping the Young diagram, e.g.
	\[
	\lambda=\begin{ytableau}
		~ & ~ & ~ \\
		~ 
	\end{ytableau}
	\qquad 
	\rightsquigarrow
	\qquad
	\lambda^*=\begin{ytableau}
		~ & ~ \\
		~ \\
		~ 
	\end{ytableau}
	\]
	
	We think of each Young diagram as the subset of the south-east quadrant in the cartesian plane, whose vertex at the north-west position sits at the origin. We identify each box with coordinates $(i,j)$ following the matrix convention, so that $i$ is the row-coordinate (from top to bottom) and $j$ is the column-coordinate (from left to right).
	For example, if $\lambda=(3,2)$, then the box in position $(1,2)$ is
	\[
	\begin{ytableau}
		~ & *(gray) ~ & ~ \\
		~ &  ~
	\end{ytableau}
	\]

	Given a Young diagram $\lambda$ and a box $(i,j)\in\lambda$, we denote by $\Hook_{i,j}(\lambda)$ the hook centered in $(i,j)$, and by $\hook_{i,j}(\lambda)$ its cardinality. For example, if $\lambda=(3,2)$, then $\Hook_{1,2}(\lambda)=\begin{ytableau}
		~ &  ~ \\
		~
	\end{ytableau}$ and $\hook_{1,2}(\lambda)=3$, while $\Hook_{2,1}(\lambda)=\begin{ytableau}
		~ &  ~ 
	\end{ytableau}$ and $\hook_{2,1}(\lambda)=2$.
	\medskip
	
	If $\lambda$ is a Young diagram, a semistandard Young tableau $T_{\lambda}$ is a filling of $\lambda$ with positive integers, weakly increasing along rows and strictly increasing down columns. For example, if $\lambda=(3,2)$, then the following are semistandard Young tableaux,
	\[  
	\begin{ytableau}
		1 & 3 & 5 \\
		2 & 4
	\end{ytableau}
	\qquad\text{ and }
	\qquad 
	\begin{ytableau}
		1 & 1 & 2 \\
		2 & 3
	\end{ytableau}.
	\]
	If its entries lie in $\{1, \dots, n\}$, its weight is $w(T_\lambda)=(w_1,\dots,w_n)$ where $w_j$ is the number of entries equal to $j$.
	In the example above, the weights are respectively $(1,1,1,1,1)$ and $(2,2,1)$.
	
	%%%%%%%%%%%%%%%%%%%%%%%%%%%%%%%%%%%%%
	\subsection{Ring of symmetric functions}
	
	We refer to \cite[Chapter~I]{Macdonald} for background on symmetric functions and their properties.
	Here we only collect the main definitions and facts needed for the rest of the paper.
	\medskip
	
	Let us fix a positive integer $n$. The ring of symmetric polynomials in $n$ variables is 
	\[ \Lambda_n=\Q[x_1,\dots,x_n]^{\mathfrak{S}_n}=\bigoplus_{d\geq0}\Lambda_n^d, \]
	where $\Lambda_n^d$ is the $\Q$-vector space of homogeneous symmetric polynomials of degree $d$ and $\mathfrak{S}_n$ acts by permuting the indeterminates $x_1, \dots, x_n$.
	
	If $\lambda$ is a Young diagram, the \emph{Schur polynomial} $\Schur_\lambda$ is the symmetric polynomial 
	\[ \Schur_\lambda(x_1,\dots,x_n)=\sum_{T_\lambda}x^{w(T_\lambda)}, \]
	where the sum runs over all semistandard Young tableau of shape $\lambda$ and entries in $\{1, \dots, n\}$ and   we write $x^{w(T_\lambda)}=x_1^{w_1}\cdots x_n^{w_n}$.
	Each Schur polynomial $\Schur_\lambda$ is homogeneous of degree $|\lambda|$. It is known that Schur polynomials of degree $d$ in $n$ variables form a basis of $\Lambda_n^d$.
	\medskip
	
	For every $m\geq n$, let us consider the natural surjective morphism
	\[ \rho_{m,n}\colon\Lambda_m\to\Lambda_n \]
	and its restriction
	\[ \rho^d_{m,n}\colon\Lambda^d_m\to\Lambda^d_n. \]
	This family of homomorphisms gives rise to an inverse system of rings, so that we can define 
	\[ \Rep^d=\varprojlim_n\Lambda^d_n. \]
	The vector space $\Rep^d$ has finite dimension equal to the number of partitions of $d$.
	
	Then, the \emph{ring of symmetric functions} is defined as 
	\[ \Rep=\bigoplus_{d\geq0}\Rep^d, \]
	and its elements are formal infinite sums of monomials. Notice that not every formal infinite sum of monomials, even if symmetric, is an element of $\Lambda$, as they may have unbounded degree (e.g.\ $x_1+x_2+x_1^2+x_2^2+x_1^3+x_2^3+\dots\notin\Lambda$).
	
	%%%%%%%%%%%%%%%%%%%%%%%%%%%%%%%%%%%%%
	\subsection{Schur and power sum functions}
	We have already defined Schur polynomials $\Schur_\lambda$ and remarked that they form an integral basis of $\Lambda_n$. 
	
	Taking the limit of the Schur polynomial $\Schur_\lambda(x_1,\dots,x_n)$ as the number of variables tends to infinity yields a well-defined symmetric function in $\Rep$, which we keep denoting by $\Schur_\lambda$. Those are called \emph{Schur functions} and they form an integral basis of $\Rep$. More precisely, $\Rep^d$ has a basis consisting of Schur functions $\Schur_\lambda$ with $|\lambda|=d$ (see \cite[I.(3.3)]{Macdonald}).
	\medskip
	
	We will need another set of generators, namely the \emph{power sum functions} 
	\[ \power_n=\sum_{i\geq1} x_i^n. \]
	By \cite[I.(2.12)]{Macdonald}, the $\power_n$ are algebraically independent over $\Q$ and 
	\[ \Rep=\Q[\power_1,\power_2,\dots]. \]
	
	\begin{remark}
		The expression of power sum functions as a combination of Schur functions is well known. Namely,
		\[ \power_n=\sum_{k=0}^{n-1} (-1)^k \Schur_{n-k,1^k}, \]
		as we shall see later in Example \ref{example:p n come schur alpha}.
	\end{remark}

	\begin{remark}\label{rmk:generators of R d}
		Each graded component $\Rep^d$ is a finite-dimensional vector space over $\Q$. Two natural bases of $\Rep^d$ are given by:
		\begin{enumerate}
			\item Schur functions $\Schur_\lambda$ with $|\lambda|=d$;
			\item monomials in power sum functions of the form $\power_\lambda$, for $|\lambda|=d$, defined as 
			\[ \power_\lambda=\power_{1}^{\blambda_1}\cdots\power_{t}^{\blambda_t}. \]
		\end{enumerate}
		Notice that Schur functions form an \emph{integral} basis, while power sum functions only form a \emph{rational} basis.
	\end{remark}
	
	%%%%%%%%%%%%%%%%%%%%%%%%%%%%%%%%%%%%%
	\subsection{Plethysms}\label{subsection:plethysms}
	
	The \emph{plethysm} is the operation on $\Rep$ analogous to composition of functions (see for example \cite{Macdonald,LO23}). In an axiomatic way, it can be defined as the unique product
	\[ \Rep\times\Rep\to\Rep,\qquad (f,g)\mapsto f[g] \]
	satisfying:
	\begin{enumerate}
		\item $\power_n[\power_k]=\power_k[\power_n]=\power_{nk}$;
		\item for every $g\in\Rep$, the assignment $f\mapsto f[g]$ is a $\Q$-algebra endomorphism of $\Rep$;
		\item the assignment $g\mapsto \power_n[g]$ is a $\Q$-algebra endomorphism of $\Rep$.
	\end{enumerate}
	Notice that by (2) we have $(f_1+f_2)[g]= f_1[g]+f_2[g]$, but in general $f[g_1+g_2]\neq f[g_1]+f[g_2]$.
	\medskip
	
	Let us list here some properties of the plethysm operation:
	\begin{itemize}
		\item if $f\in\Rep^d$ and $g\in\Rep^{e}$, then $f[g]\in\Rep^{de}$;
		\item for every $f\in\Rep$, we have $\power_n[f]=f[\power_n]$;
		\item (associativity) for any $f,g,h\in\Rep$, we have $f[g[h]]=(f[g])[h]$ (\cite[I.(8.7)]{Macdonald});
		\item if $\Schur_\lambda[\Schur_\mu]=\sum_{\nu}c^\nu_{\lambda, \mu} \Schur_\nu$, then $c^\nu_{\lambda, \mu}\geq0$ (\cite[Appendix~I.A]{Macdonald}).
	\end{itemize}

	%%%%%%%%%%%%%%%%%%%%%%%%%%%%%%%%%%%%%
	\subsection{Hall product}\label{section:Hall product}
	It is possible to define a symmetric and positive definite bilinear form on $\Rep$
	\[ \langle\cdot,\cdot\rangle\colon\Rep\times\Rep\longrightarrow\Q, \]
	uniquely determined by the property that 
	\[ \langle\Schur_\lambda,\Schur_\mu\rangle=\delta_{\lambda,\mu}, \]
	where $\delta_{\lambda,\mu}$ is the Kronecker delta function. In other words, Schur functions form an orthonormal basis of $(\Rep,\langle\cdot,\cdot\rangle)$ (see \cite[Section~I.4]{Macdonald}). By bilinearity, elements of $\Lambda^d$ are orthogonal to those of $\Lambda^e$ as soon as $d \neq e$.
	
	Monomials of power sums are also pairwise orthogonal, but not orthonormal, with respect to $\langle\cdot,\cdot\rangle$.
	\begin{lemma}[\protect{\cite[I.(4.7)]{Macdonald}}]\label{lemma:Hall products of monomials in power sums}
		We have
		\[ \langle \power_\lambda,\power_\mu\rangle=\left\{\begin{array}{ll}
			z_\lambda:=1^{\blambda_1}\cdots t^{\blambda_t}\cdot\blambda_1!\cdots\blambda_t! & \text{if } \lambda=\mu \\
			0 & \text{if } \lambda\neq\mu\, .
		\end{array}\right.\] \qed
	\end{lemma}

	%%%%%%%%%%%%%%%%%%%%%%%%%%%%%%%%%%%%%
	%%%%%%%%%%%%%%%%%%%%%%%%%%%%%%%%%%%%%
	\section{Multiplication and differentiation by power sum functions}\label{section:derivation}
	
	If $\varphi \in \Rep$, let $\varphi^\perp\colon \Rep \to \Rep$ denote the $\Q$-linear operator adjoint to multiplication by $\varphi$. More precisely, $\varphi^\perp$ is the unique endomorphism of $\Rep$ satisfying
	\[ \langle \varphi^\perp(f), g\rangle = \langle f, \varphi g\rangle \]
	for all $f, g \in \Rep$. We will primarily be interested in the special case when $\varphi=\power_n$. The following lemma recalls the well-known fact that $\power_n^\perp$ is a derivation of $\Rep$.
	
	\begin{lemma}\label{lemma:D aggiunto P}
		We have 
		\[ \power_n^\perp=n\frac{\partial}{\partial\power_n}\,. \]
		\begin{proof}
			By the linearity of both $n\frac{\partial}{\partial\power_n}$ and the Hall inner product, it suffices to show that 
			\[ \langle n\frac{\partial}{\partial\power_n}\power_\lambda,\power_\mu\rangle=\langle \power_\lambda,\power_n\power_\mu\rangle. \]
			For a partition $\lambda$, we denote by $\lambda^{(n)}$ and $\lambda_{(n)}$ the partitions such that, respectively,
			\[ \blambda^{(n)}_k=\begin{cases}
				\blambda_k & \mbox{if } k\neq n \\
				\blambda_n+1 & \mbox{if } k=n
			\end{cases}
			\qquad\mbox{and}\qquad
			(\blambda_{(n)})_k=\begin{cases}
				\blambda_k & \mbox{if } k\neq n \\
				\blambda_n-1 & \mbox{if } k=n
			\end{cases}
			\]
			
			It is straightforward to see that $n\frac{\partial}{\partial\power_n}\power_\lambda=n\blambda_n\power_{\lambda_{(n)}}$, while $\power_n(\power_\mu)=\power_{\mu^{(n)}}$. Therefore,
			\[ \langle n\frac{\partial}{\partial\power_n}\power_\lambda,\power_\mu\rangle=n\blambda_n\langle\power_{\lambda_{(n)}},\power_\mu\rangle\qquad\text{and}\qquad \langle \power_\lambda,\power_n \power_\mu\rangle=\langle\power_\lambda,\power_{\mu^{(n)}}\rangle, \]
			and the equality follows immediately from Lemma \ref{lemma:Hall products of monomials in power sums}.
		\end{proof}
	\end{lemma}
	
	As previously noted, the operators of multiplication and differentiation by $\power_n$ are ubiquitous in the theory of symmetric functions and have been extensively studied. A graphical description of their action on Schur functions is well known. For instance, multiplication by $\power_n$ is described in \cite[Formula (2), p.~48]{Macdonald} as
	\begin{equation}\label{eqn:Murnaghan Nakayama} \power_n\Schur_\lambda=\sum_\mu (-1)^{\operatorname{ht}(\mu\setminus\lambda)}\Schur_\mu, \end{equation}
	where the sum runs over all Young diagrams $\mu$ containing $\lambda$ such that $\mu\setminus\lambda$ is a border strip of size $n$. Here, $\operatorname{ht}(\mu\setminus\lambda)$ denotes the height of the strip.
	
	By adjointness, we obtain a similar description for $\power_n^\perp$, where instead of \emph{adding} border strips, one \emph{removes} them or, equivalently, collapses hooks. Since we require this description, we conclude this section by providing a proof. We start with a working definition.
	
	\begin{definition}\label{defn:alpha i j}
		Let $\lambda$ be a Young diagram and let $(i,j)\in\lambda$ be a box. Then 
		\[ \lambda_{(i,j)}:=\lambda\setminus\Hook_{i,j}(\lambda) \]
		is obtained from $\lambda$ by collapsing the hook centered at position $(i,j)$. If the hook disconnects $\lambda$, we shift the lower component up and to the left. 
	\end{definition}
	
	\begin{example}
		Consider $\lambda=(3,2,2)=\begin{ytableau}
			~ & ~ & ~  \\
			~ &  ~  \\
			~ & ~
		\end{ytableau}$ . Then 
		\[
		\lambda_{(1,2)}=\begin{ytableau}
			~   \\
			~  \\
			~
		\end{ytableau}
		\qquad\text{ and }
		\qquad
		\lambda_{(2,1)}=\begin{ytableau}
			~ & ~ & ~  \\
			~  
		\end{ytableau}
		\]
	\end{example}
	
	\begin{proposition}[Geometric description of $\power_n^\perp$]\label{prop:geometric interpretation of D}
		Let $\operatorname{H}_n(\lambda)$ denote the set of boxes $(i,j)\in\lambda$ such that $\hook_{i,j}(\lambda)=n$. Then 
		\begin{equation*}    \power_n^\perp(\Schur_\lambda)=\sum_{(i,j)\in\operatorname{H}_n(\lambda)} (-1)^{\varepsilon(i,j)} \Schur_{\lambda_{(i,j)}}, 
		\end{equation*}
		where $\varepsilon(i,j)$ is the leg-length, i.e., the number of boxes sitting strictly below $(i,j)$.
		\begin{proof}
			First, by adjointness, $\Schur_\mu$ appears as a summand in $\power_n^\perp(\Schur_\lambda)$ if and only if $\Schur_\lambda$ appears as a summand in $\power_n(\Schur_\mu)$. By \eqref{eqn:Murnaghan Nakayama}, the latter condition holds if and only if $\lambda$ contains $\mu$ and $\lambda/\mu$ is a border strip of size $n$. Crucially, this is equivalent to saying that $\mu=\lambda_{(i,j)}$ for some $(i,j)\in\operatorname{H}_n(\lambda)$. This is clear when the border strip is itself a hook; if it is not, the associated hook disconnects the Young diagram, and the claim still holds since we agreed to shift the lower component up and to the left. 
			
			Finally, the sign conventions are easily seen to agree, which completes the proof.
		\end{proof}
	\end{proposition}
	
	%%%%%%%%%%%%%%%%%%%%%%%%%%%%%%%%%%%%%
	\subsection{Connection with the representation theory of symmetric groups}\label{section:symmetric group I}
	
	If $\lambda$ is a partition, set $z_\lambda = \prod_{k \ge 1} \blambda_k! k^{\blambda_k}$. When $|\lambda|=n$, then $z_\lambda$ is the order of the centralizer in the symmetric group $\Ss_n$ of a permutation of cycle type $\lambda$. Equivalently, $\Ss_n$ contains $n!/z_\lambda$ elements of cycle type $\lambda$. We denote the corresponding conjugacy class by $C_\lambda$.
	
	The space $\Cl(\Ss_n)$ of $\Q$-valued class functions on the symmetric group $\Ss_n$ is an inner product space with respect to
	\[ \langle \phi,\psi\rangle = \frac{1}{n!} \sum_{x \in \Ss_n} {\phi(x)}\psi(x). \]
	The irreducible (and $\Z$-valued) characters $\chi^\lambda$ for a partition $\lambda$ with $|\lambda|=n$ form an orthonormal basis for $\Cl(\Ss_n)$.
	If $g\in \Ss_n$ has cycle type $\mu$, we set $\chi^\lambda(\mu):= \chi^\lambda(g)$, following the notation of \cite[\S 7.18]{StanleyII}.
	The \emph{Frobenius characteristic map} $\ch_n$ from the space $\Cl(\Ss_n)$ of class functions on the symmetric group $\Ss_n$ to $\Lambda^n$ is defined by linearly extending
	\[ [C_\mu] \stackrel{\ch_n}{\longmapsto} \frac{\power_\mu}{z_\mu}, \]
	where $[X]$ denotes the characteristic function of any given subset $X \subset \Ss_n$. 
	We observe that
	\[ \Bigl\langle \frac{\power_\mu}{z_\mu}, \frac{\power_\mu}{z_\mu} \Bigr\rangle = \frac{1}{z_\mu} \]
	and correspondingly,
	\[ \bigl\langle \, [C_\mu], [C_\mu] \, \bigr\rangle = \frac{|C_\mu|}{n!} = \frac{1}{z_\mu}, \]
	as there are $n!/z_\mu$ elements in $\Ss_n$ with cycle type $\mu$.
	
	Since the characteristic functions of distinct conjugacy classes are clearly orthogonal, the map $\ch_n$ is an isometry. Furthermore, by \cite[Theorem~7.18.5]{StanleyII} or \cite[p.~114, (7.5)]{Macdonald},
	\[ \chi^\lambda = \sum_{|\mu|= n} \chi^\lambda(\mu) [C_\mu] \mapsto \sum_{|\mu|=n} \chi^\lambda(\mu)\frac{\power_\mu}{z_\mu} = \Schur_\lambda. \]
	For example, $\Schur_n = \Schur_{(n)} = \sum_\mu \power_\mu/z_\mu$ corresponds to the trivial character of $\Ss_n$.
	We can now translate formula \eqref{eqn:Murnaghan Nakayama} in terms of characters of the symmetric group. We obtain that for any $g \in \Ss_{n-k}$,
	\begin{equation} \chi^\mu \bigl( g \cdot (n-k+1,\ldots, n)\bigr)=\sum_\lambda (-1)^{\operatorname{ht}(\mu/\lambda)} \chi^\lambda(g), \label{eq:MN} \end{equation}
	where $(n-k+1,\ldots, n)$ is the cyclic permutation on the last $k$ letters and, as before, the sum runs over all Young diagrams $\mu$ containing $\lambda$ such that $\mu/\lambda$ is a border strip with $k$ boxes. This is the celebrated Murnaghan--Nakayama rule for computing the character table of symmetric groups.

	%%%%%%%%%%%%%%%%%%%%%%%%%%%%%%%%%%%%%
	%%%%%%%%%%%%%%%%%%%%%%%%%%%%%%%%%%%%%
	\section{A remarkable derivation of the algebra of symmetric functions}\label{section:remarkable derivation}
	
	Recall that each homogeneous component $\Rep^d$ is spanned by the power sum monomials $\power_\lambda$ satisfying $|\lambda|=d$. Since only finitely many of the partial derivatives $\frac{\partial}{\partial\power_n}$ act non-trivially on any given $\power_\lambda$, we have a well-defined (non-homogeneous) derivation
	\begin{equation}\label{defn:D}
		\D\colon\Rep\to\Rep,\qquad\D:=\sum_{n\geq1}\power_n^\perp=\sum_{n\geq1}n\frac{\partial}{\partial\power_n}.
	\end{equation}
	
	The operator $\D$ can also be thought of as the unique derivation on $\Rep \cong \Q[\power_1, \power_2, \dots]$ satisfying $\D(\power_n) = n$ for all $n \in \N$.
	
	\begin{remark}\label{rmk:D gia nota}
		The operator $\D$ is a natural object of study and has already appeared in the literature in various guises. For example, in \cite[Section~2.2]{Yanagida}, the author defines an operator $D_q(z)$ that reduces to our $\D$ after dividing by $1-q$ and taking the limit as $q\to1$. 
	\end{remark}
	
	Our first result is the following compatibility of the operator $\D$ with the Hall inner product.
	
	\begin{thm}\label{thm:quasi isometria}
		For every $d>0$, the derivation $\D$ is a quasi-isometry on $\Rep^d$; that is, for all homogeneous symmetric functions $f,g\in\Rep^d$, 
		\[ \langle\D(f),\D(g)\rangle=d\langle f,g\rangle. \]
		In particular, $\D$ acts injectively on $\Rep^d$ for every $d>0$.
		\begin{proof}
			First, we claim that 
			\[ \langle\D(f),\D(g)\rangle=\sum_{n\geq1}\langle\power_n^\perp(f),\power_n^\perp(g)\rangle. \]
			Indeed, writing $\power_n^\perp(\power_\lambda)=n\blambda_n\power_{\lambda_{(n)}}$ as in the proof of Lemma \ref{lemma:D aggiunto P}, by linearity and Lemma \ref{lemma:Hall products of monomials in power sums}, we obtain 
			\[ \langle\power_n^\perp(f),\power_m^\perp(g)\rangle=0\qquad\mbox{ if } n\neq m\,. \]
			From the definition of $\power_n^\perp$ as the adjoint to multiplication by $\power_n$, we then obtain
			\[ \langle\D(f),\D(g)\rangle=\sum_{n\geq1}\langle\power_n^\perp(f),\power_n^\perp(g)\rangle=\langle f,\sum_{n\geq1}\power_n(\power_n^\perp(g))\rangle = \langle f, \sum_{n \geq 1} n\,\power_n \frac{\partial g}{\partial \power_n}\rangle, \]
			and the claim follows immediately upon noticing that
			$$\sum_{n \geq 1} n \, \power_n \frac{\partial}{\partial \power_n}$$
			is the grading operator on $\Rep$ that maps each $g\in \Rep^d$ to $dg$.
		\end{proof}
	\end{thm}
	
	\begin{remark}
		Although $\D$ is injective on each $\Rep^d$ for $d>0$, it fails to be injective on the ideal $\Rep^{>0}$. For example, $\D(\power_n/n) = 1$ for every $n>0$. 
		
		On the other hand, $\D\colon\Rep\to\Rep$ is surjective. The following argument is due to the referee and simplifies our previous (inductive) argument. First of all, notice that for every $n\geq1$ we have $\D(\power_n-n\power_1)=0$. Performing the change of variables $\power_n'=\power_n-n\power_1$ for every $n\geq2$, we have $\Rep\cong\Q[\power_1,\power_2',\cdots]$ and every $f\in\Rep$ can be written as 
		\[ f=\sum_{j=0}^k a_j(\power_2',\power_3',\cdots)\power_1^j \]
		where $a_j(\power_2',\power_3',\cdots)\in\Q[\power_2',\power_3',\cdots]$. Since $\D(\power_1)=1$ and $\D(\power_n')=0$ for $n\geq2$, it follows that 
		\[ f=\D\left(\sum_{j=0}^k a_j(\power_2',\power_3',\cdots)\frac{\power_1^{j+1}}{j+1} \right). \]
	\end{remark}
	
	\begin{remark}\label{rmk:inverse of D}    
		An explicit inverse to the restriction $\D|_{\Rep^d}$ is easily constructed. Indeed, since
		$$\sum_{n\geq 1}n\power_n \frac{\partial}{\partial \power_n} = \sum_{n \geq 1} \power_n \power_n^\perp$$
		acts via multiplication by $d$ on $\Rep^d$, we have 
		$$f = \frac{1}{d} \sum_n \power_n^\perp(f)\, \power_n .$$
	\end{remark}
	
	Applying Proposition \ref{prop:geometric interpretation of D} provides a useful method for computing $\D$ on Schur functions. Explicitly, we have
	\begin{equation}\label{eqn:geometric expression of D}
		\D(\Schur_\lambda)=\sum_{(i,j)\in\lambda} (-1)^{\varepsilon(i,j)} \Schur_{\lambda_{(i,j)}}, 
	\end{equation}
	where $\varepsilon(i,j)$ is the number of boxes sitting strictly below $(i,j)$.
	
	\begin{example}\label{example:several examples of D}
		We collect several examples below:
		\begin{enumerate}
			\item $\D(\Schur_m)=\Schur_{m-1}+\Schur_{m-2}+\cdots+1$;
			\item $\D(\Schur_{1^m})=\Schur_{1^{m-1}}-\Schur_{1^{m-2}}+\cdots+(-1)^{m-1}1$;
			\item $\D(\Schur_{2,1})=\Schur_2+\Schur_{1,1}-1$;
			\item $\D(\Schur_{3,2,2})=\Schur_{3,2,1}+\Schur_{2,2,2}+\Schur_{3,2}-\Schur_{3,1,1}-\Schur_{3,1}+\Schur_{1,1,1}+\Schur_{1,1}$.
		\end{enumerate}
	\end{example}
	
	\begin{example}
		By Example \ref{example:several examples of D} and Remark \ref{rmk:inverse of D}, it follows that
		\[ \Schur_m=\frac{1}{m}\sum_{k=1}^{m} \Schur_{m-k}\power_{k} \qquad \mbox{ and }\qquad \Schur_{1^m}=\frac{1}{m}\sum_{k=1}^{m} (-1)^{k-1}\Schur_{1^{m-k}}\power_{k}  \]
		(cf.\ \cite[Formulas (2.11) and (2.11') in Section~I]{Macdonald}).
	\end{example}
	
	\begin{example}\label{example:p n come schur alpha}
		We can use Lemma \ref{lemma:D aggiunto P} and Proposition \ref{prop:geometric interpretation of D} to express power sum functions in terms of Schur functions. Writing $\power_n=\sum_\lambda m_\lambda\Schur_\lambda$, we have 
		\[ m_\lambda=\langle\power_n,\Schur_\lambda\rangle=\langle 1, \power_n^\perp \Schur_\lambda \rangle.
		\]
		The latter is non-zero only if $\lambda$ is a hook partition of size $n$, in which case $\power_n^\perp(\Schur_\lambda)=\pm 1$. Taking care of the signs yields
		\[ \power_n=\sum_{k=0}^{n-1} (-1)^k \Schur_{n-k,1^k}. \]
	\end{example}
	
	\begin{example}\label{example:D of transpose}
		Let $\lambda$ be a partition and $\lambda^*$ its conjugate. 
		We have 
		\[ \D(\Schur_{\lambda^*})=\sum_{(i,j)\in\lambda}(-1)^{\varepsilon(i,j)+\eta(i,j)}\, \Schur_{(\lambda_{(i,j)})^*}, \]
		where $\eta(i,j)$ is $1$ if $\hook_{ij}(\lambda)$ is even and $0$ if $\hook_{ij}(\lambda)$ is odd (i.e., there is a sign change only for even values of the hook number).
		
		For example, if $\lambda=(3,1)$, then
		\[ \D(\Schur_{3,1})=\Schur_3+\Schur_{2,1}+\Schur_{1,1}-1\quad\text{ and }\quad \D(\Schur_{2,1,1})=\Schur_{1,1,1}+\Schur_{2,1}-\Schur_{2}+1. \]
	\end{example}
	
	%%%%%%%%%%%%%%%%%%%%%%%%%%%%%%%%%%
	\subsection{A recursion for Littlewood--Richardson coefficients}
	
	Since the operator $\D$ is a derivation on $\Rep$, as a corollary of Theorem \ref{thm:quasi isometria}, we obtain a recursive formula to express the product of two Schur functions as a sum of Schur functions. The resulting recursion is, however, more involved than what might be obtained by using the adjointness of a single $n\partial/\partial \power_n$ with the corresponding multiplication by $\power_n$.
	
	\begin{corollary}\label{cor:LR}
		Let $\lambda$ and $\mu$ be two partitions. If we write $\Schur_\lambda\Schur_\mu=\sum_{\nu}m_{\lambda,\mu}^\nu\Schur_\nu$, then 
		\[ m_{\lambda,\mu}^\nu=\frac{1}{d}\langle\D(\Schur_\lambda)\Schur_\mu+\Schur_\lambda\D(\Schur_\mu),\D(\Schur_\nu)\rangle, \]
		where $d=|\lambda|+|\mu|$.
		\begin{proof}
			The proof follows at once from Theorem \ref{thm:quasi isometria} and the fact that $\D$ is a derivation.
		\end{proof}
	\end{corollary}
	
	\begin{remark}
		The right-hand side of the formula above requires the computation of several products. Nevertheless, these are products of Schur functions of degree strictly smaller than $d$, meaning they are known recursively.
		
		Recall that the coefficients $m_{\lambda,\mu}^\nu$ can be computed in a combinatorial and direct way using the Littlewood--Richardson rule (e.g., \cite[Section~I.9]{Macdonald}). 
	\end{remark}
	
	\begin{example}
		Let us compute the product $\Schur_2\Schur_1$. First, recall that $\D(\Schur_2)=\Schur_1+1$ and $\D(\Schur_1)=1$. Therefore,
		\[ \D(\Schur_2\Schur_1)=(\Schur_1+1)\Schur_1+\Schur_2=2\Schur_2+\Schur_{1,1}+\Schur_1, \]
		where we have used the well-known decomposition $\Schur_1^2=\Schur_2+\Schur_{1,1}$.
		Now,
		\begin{align*}
			\D(\Schur_3) = \Schur_2+\Schur_1+1 \Longrightarrow & \; \langle \Schur_2\Schur_1,\Schur_3\rangle = \frac{1}{3}(2+1)=1 \\
			\D(\Schur_{2,1}) = \Schur_2+\Schur_{1,1}-1 \Longrightarrow & \; \langle \Schur_2\Schur_1,\Schur_{2,1}\rangle = \frac{1}{3}(2+1)=1 \\
			\D(\Schur_{1,1,1}) = \Schur_{1,1}-\Schur_1 +1 \Longrightarrow & \; \langle \Schur_2\Schur_1,\Schur_{1,1,1}\rangle = \frac{1}{3}(1-1)=0
		\end{align*}
		from which it follows that 
		\[ \Schur_2\Schur_1=\Schur_3+\Schur_{2,1}. \]
	\end{example}
	
	%%%%%%%%%%%%%%%%%%%%%%%%%%%%%%%%%%%%%%%%
	\subsection{Action of $\D$ and class functions}\label{section:symmetric group II}
	
	With the same notation as in Section \ref{section:symmetric group I}, we have
	\[ \frac{\power_\lambda}{z_\lambda} \stackrel{\D}{\longmapsto} \D\Bigl( \frac{\power_\lambda}{z_\lambda}
	\Bigr)
	= \frac{1}{z_\lambda} \sum_{k \ge 1} k \,\blambda_k \,\power_{\lambda_{(k)}} = \sum_{k \ge 1}
	\frac{\power_{\lambda_{(k)}}}{z_{\lambda_{(k)}}} \]
	where $\lambda_{(k)}$ is as in the proof of Lemma \ref{lemma:D aggiunto P}, and we agree to set $\power_{\lambda_{(k)}} = 0$ if $\blambda_k = 0$.
	
	Conjugating $\D$ by $\ch$ yields the unique linear map $\D_\Cl \colon \Cl(\Ss_n) \to \bigoplus_{k \ge 1} \Cl(\Ss_{n-k}) $ satisfying
	\begin{equation}
		\label{eq:DClFirst}
		[C_\lambda] \stackrel{\D_\Cl}{\longmapsto} \sum_{k \ge 1} [C_{\lambda_{(k)}}],
	\end{equation}
	where, once again, $[C_{\lambda_{(k)}}] = 0$ if $\blambda_k = 0$.
	
	Let $(n-k+1,\ldots, n) \in \Ss_n$ denote the $k$-cycle permuting the last $k$ letters.
	We can rewrite \eqref{eq:DClFirst} as
	\[ [C_\lambda] \stackrel{\D_\Cl}{\longmapsto}\sum_{k \ge 1} \left[\{ g \in \Ss_{n-k} \,|\, g(n-k+1,\ldots, n) \in C_\lambda\}\right].\]
	
	Hence, by linearity, for any $\phi \in \Cl(\Ss_n)$, we have
	\[\D_\Cl (\phi) = \sum_{k=1}^n \Bigl( g \in \Ss_{n-k} \mapsto \phi\bigl( g (n-k+1, \ldots, n) \bigr) \Bigr), \]
	where we use the function notation $g \mapsto f(g)$ to denote $f$.
	
	This shows that $\D_\Cl$ can be interpreted as a type of restriction map: the component of the image of $\phi$ in $\Cl(\Ss_{n-k})$ is the function induced by $\phi$ on elements of $\Ss_{n-k}$ via the embedding $\Ss_{n-k} \to \Ss_n$ defined by $g \mapsto g (n-k+1,\ldots, n)$.
	Writing $\phi|_{\Ss_{n-k} (n-k+1,\ldots, n)}$ for this restricted function, we have
	\begin{equation} \label{eq:DCl} \D_\Cl(\phi) = \sum_{k \ge 1} \phi|_{\Ss_{n-k} (n-k+1,\ldots, n)}, \end{equation}
	where the terms on the right-hand side can be computed using the Murnaghan--Nakayama rule in \eqref{eq:MN}.
	Since any class function on $\Ss_n$ is uniquely determined by its values on elements of the form $g (n-k+1,\ldots, n)$ for $1 \le k \le n$, it follows that $\D_\Cl(\phi)$ determines $\phi$. This establishes the injectivity of $\D_\Cl$ on $\Cl(\Ss_n)$, and correspondingly, the injectivity of $\D$ when restricted to $\Rep^n$.
	
	Similar constructions can be found in the literature. For instance, the recurrence for the plethysm coefficients $c^\nu_{(n),(m)}$ in \cite[Proposition 5.1]{EPW} is proved by conditioning on the cycle containing $n$ in a permutation in $\Ss_n$, thus mirroring \eqref{eq:DCl}.
	
	%%%%%%%%%%%%%%%%%%%%%%%%%%%%%%%%%%%%%
	%%%%%%%%%%%%%%%%%%%%%%%%%%%%%%%%%%%%%
	\section{A recursive formula for plethysms}\label{section:plethysms}
	
	Multiplication operators $m_g: f \mapsto fg = gf$ by any given element $g \in \Rep$ are $\Q$-linear endomorphisms of $\Rep$; similarly, (left) plethysm maps $f \mapsto f[g]$ lie in $\End_\Q \Rep$ for all choices of $g \in \Rep$.
	In particular, all {\em Adams operations} 
	\[ \Psi^n\colon f \mapsto \power_n[f] = f[\power_n] \] 
	lie in $\End_\Q \Rep$ as they are indeed $\Q$-algebra homomorphisms. One has the following identities:
	$$
	m_f \circ m_g = m_{fg}, \qquad \Psi^h \circ \Psi^k = \Psi^{hk}, \qquad \Psi^h \circ m_f = m_{\Psi^h(f)},
	$$
	for all $f, g \in \Rep$ and $h, k \in \N$. In all that follows we shall abuse notation and denote $m_f$ simply by $f$ thus implicitly identifying $\Rep$ with a subring of $\End_\Q \Rep$, so that
	\begin{equation}\label{defn:prodotto strano}
		(f \circ \Psi^h) \circ (g \circ \Psi^k) = (f\, \Psi^h(g)) \circ \Psi^{hk}.
	\end{equation}
	Notice that both $1 \in \Rep$ and $\Psi^1$ are the identity endomorphism, so that $1 \circ \Psi^h = \Psi^h$ and $f \circ \Psi^1 = f$ for all $h \in \N, f \in \Rep$.
	
	We are going to modify the definition of $\D$ to a map $\DD$ taking values in $\End_\Q \Rep$, whose compatibility with the plethysm operation will be the main focus in the rest of this section.
	
	\begin{definition}\label{defn:Ds}
		Define $\DD\colon\Rep\to\End_\Q \Rep$ as
		\[ \DD(f)=\sum_{n\geq1} \power_n^\perp(f) \circ \Psi^n\,. \]
	\end{definition}

	%\alessandro{Non ho scritto la dimostrazione del fatto che la somma nel Remark 4.2 è diretta. La scrivo qui e mi dite se va inclusa (non penso). Supponete di avere una relazione lineare FINITA $\sum f_n \circ \Psi^n = 0$ e sia $N$ il più grande intero tale che $\power_N$ NON COMPAIA in nessun $f_n$ quando vengono espressi come polinomi nei powersum. Ora scegliete $k$ in modo che $kd>N$ dove $d$ è il massimo $n$ tale che $f_n \neq 0$. Applicando la relazione lineare a $\power_k$, l'ultimo addendo non nullo è $f_d \power_{kd}$ che non può cancellarsi con nessun altro termine.}
	
	\begin{remark}
		A routine check makes sure that the sum $\Rep \Psi: = \sum_{n \geq 1} \Rep \circ \Psi^n$ is direct and \eqref{defn:prodotto strano} shows that $\Rep\Psi \subset \End_\Q \Rep$ is a subring with respect to composition. Notice that $\DD: \Rep \to \Rep \Psi$ is homogeneous of degree $0$ as soon as we set $\deg \Psi^n = n$.
		
		$\Rep$-linear independence of the $\Psi^n$ provides a more abstract interpretation of $(\Rep\Psi, \circ)$ as follows. Denote by $\Psi = \{\Psi^n, n \geq 1\}\simeq \N$ the semigroup of all Adams operations under composition. The natural action of $\Psi$ on the ring $\Rep$ extends to a (Hopf) action of the semigroup (bi)algebra $\Q\Psi$, where the bialgebra structure is obtained by imposing all Adams operations to be group-like. The corresponding smash product $\Rep \,\sharp\, \Q\Psi$ coincides, as a ring, with the $\Q$-algebra $\Rep\Psi$ endowed with the above composition product $\circ$.
	\end{remark}
	
	\begin{remark}\label{twoproducts}
		We stress the fact that we have another different and naturally defined product $\cdot$ on the vector space $\Rep\Psi$, obtained by $\Rep$-linearly extending composition on $\Psi$. Then $(\Rep\Psi, \cdot)$ coincides with the semigroup algebra of the semigroup $\Psi$ with coefficients in the commutative ring $\Rep$ and $\DD = \sum_{n \geq 1} \Psi^n \power_n^\perp$ is a derivation of $(\Rep\Psi,\cdot)$ as $\DD$ is a locally finite $\Rep\Psi$-linear combination of the derivations $\power_n^\perp$, and $(\Rep\Psi, \cdot)$ is commutative. This observation will be of paramount importance in the proof of 
		Theorem \ref{thm:plethysm product} below.
	\end{remark}
	By Proposition \ref{prop:geometric interpretation of D} we can combinatorially compute $\DD$ of a Schur function,
	\begin{equation}\label{eqn:D s}
		\DD(\Schur_\lambda)=\sum_{(i,j)\in\lambda} (-1)^{\epsilon(i,j)} \Schur_{\lambda_{(i,j)}}\, \circ \Psi^{\hook_{i,j}(\lambda)}, 
	\end{equation}
	where $\varepsilon(i,j)$ is the number of boxes strictly below $(i,j)$.
	\medskip
	
	We can define a generalised Hall product on $\Rep\Psi$ by bilinearly extending
	\[ \llangle f\circ \Psi^h, g\circ \Psi^k\rrangle:=\langle f,g\rangle \,\delta_{h,k}, \]
	where $\delta_{h,k}$ is the Kronecker delta.
	The following is a verbatim reformulation of Theorem \ref{thm:quasi isometria}.
	
	\begin{thm}\label{thm:quasi isometria s}
		For any $f,g\in\Rep^d$, we have
		\[ \llangle \DD(f),\DD(g)\rrangle=d\langle f,g \rangle.  \] 
		Moreover, the restriction of $\DD$ to each $\Rep^d, d>0,$  is injective. \qed
	\end{thm}
	
	We also partially extend the plethysm operation to $\Rep\Psi$ in the following way.
	If $f=f_1\circ \Psi^1+\cdots+f_n\circ \Psi^n\in\Rep\Psi$ and $g\in\Rep$, then we set
	\[ f[g]:=f_1[g]\circ \Psi^1+\cdots+f_n[g]\circ \Psi^n \qquad\text{and}\qquad g[f]:=g[f_1]\circ \Psi^1+\cdots + g[f_n]\circ \Psi^n. \]
	We will only apply these operations in cases where they are defined.

	The main result of this section is the following result, which is analogous to the chain rule for derivatives.
	
	\begin{thm}[Chain rule for plethysms]\label{thm:plethysm product}
		For any $f,g\in\Rep$ we have
		\begin{equation}\label{eqn:D di un pletismo}
			\DD(f[g])=(\DD(f)[g])\circ \DD(g). 
		\end{equation}
		\begin{proof}
			First of all, notice that the claim follows when $f=\power_h$ and $g=\power_k$. In fact, as recalled in Section \ref{subsection:plethysms}, we have $\power_h[\power_k]=\power_{hk}$ and it is easy to see that both sides of \eqref{eqn:D di un pletismo} compute to $hk\circ \Psi^{hk}$.
			
			Set $f=\power_h$ and consider the ring endomorphism $\Psi^h: g \mapsto \power_h[g]$. In Remark \ref{twoproducts} we noticed that $\DD$ is a derivation on $(\Rep\Psi,\cdot)$; it follows that the application $g \mapsto \DD(\power_h[g])$ is a $\Psi^h$-\emph{twisted derivation}, i.e.\ 
			\[ \DD(\Psi^h(g_1g_2))=\Psi^h(g_1)\cdot\DD(\Psi^h(g_2))+\DD(\Psi^h(g_1))\cdot\Psi^h(g_2). \]
		
			As for the right-hand side, $\DD(\power_h) = h\circ \Psi^h$ so that also $\DD(\power_h)[g] = h[g]\circ \Psi^h = h \circ \Psi^h$. Then $(\DD(\power_h)[g]) \circ \DD(g) = h\power_h[\DD(g)]$ is again a $\Psi^h$-twisted derivation $\Rep \to (\Rep\Psi, \cdot)$.
			
			Since $\Rep=\Q[\power_1,\power_2,\dots]$ is generated as a $\Q$-algebra by elements $\power_k$, and both sides are $\Psi^h$-twisted derivations $\Rep \to (\Rep\Psi, \cdot)$ that agree when evaluated at $g = \power_k$ for every $k \in \N$, we conclude that 
			\[ \DD(\power_h[g])=(\DD(\power_h)[g])\circ \DD(g),  \]
			for all $g\in \Rep$.
			
			Finally, for every fixed $g\in\Rep$ let us consider the ring endomorphism $[g]: f\mapsto f[g]$. Again by Remark \ref{twoproducts}, we see that  both sides of \eqref{eqn:D di un pletismo} are $[g]$-\emph{twisted derivations}.
			%, i.e.\[ \DD(f_1f_2[g])=f_1[g]\cdot\DD(f_2[g])+\DD(f_1[g])\cdot f_2[g]. \]
			%Finally, for every fixed value of $g\in\Rep$,  both sides of (\ref{eqn:D di un pletismo}) are derivations $\Rep \to (\Rep\Psi, \cdot)$, as functions of $f$, by distributivity of $\circ$ and the fact that left plethysms are ring endomorphisms.
			%In fact, the right-hand side is clearly a derivation, while for the left-hand side we can again assume that $f=\power_\lambda$ and use that the assignment $f\mapsto f[g]$ is an endomorphism of $\Rep$ to conclude (see \Cref{subsection:plethysms} again). 
			%Moreover, also the right-hand side of (\ref{eqn:D di un pletismo}) is a $[g]$-twisted derivation.
			As before, the two sides coincide when $f = \power_h$ for every $h \in \N$, and this ensures that \eqref{eqn:D di un pletismo} holds for all choices of $f \in \Rep$, thus finishing the proof.
		\end{proof}
	\end{thm}
	As a straightforward consequence, we get a recursive formula for plethysms.
	
	\begin{corollary}\label{cor:plethysm formula}
		If $\lambda$ and $\mu$ are two partitions and we write $\Schur_\lambda[\Schur_\mu]=\sum_{\nu} c^\nu_{\lambda,\mu}\,\Schur_\nu$, then 
		\[ c^\nu_{\lambda,\mu}=\frac{1}{d}\llangle \DD(\Schur_\lambda)[\Schur_\mu])\circ\DD(\Schur_\mu),\DD(\Schur_\nu)\rrangle, \] 
		where $d=|\lambda|\,|\mu|$. \qed
	\end{corollary}
	
	Let us see it explicitly in some examples.
	
	\begin{example}\label{example:S2 plet S2 seconda parte}
		Let us compute the plethysm $\Schur_2[\Schur_2]$. First of all, $\DD(\Schur_2)=\Schur_1\circ \Psi^1+\Psi^2$. Then
		\[ \DD(\Schur_2)[\Schur_2]= \Schur_2\circ \Psi^1+\Psi^2 \]
		so that 
		\begin{align*}
			\DD(\Schur_2[\Schur_2]) & =\, (\Schur_2\circ \Psi^1+\Psi^2)\circ(\Schur_1\circ \Psi^1+\Psi^2) \\
			& = \, \Schur_2\Schur_1\circ \Psi^1+(\Schur_2+\power_2)\circ \Psi^2+\Psi^4 \\
			& = \, (\Schur_3+\Schur_{2,1})\circ \Psi^1 + (2\Schur_2-\Schur_{1,1})\circ \Psi^2 + \Psi^4\, .
		\end{align*}
		Explicitly, we can compute
		\begin{align*}
			\DD(\Schur_4) & = \, \Schur_3\circ \Psi^1+\Schur_2\circ \Psi^2+\Schur_1\circ \Psi^3 + \Psi^4 \\
			\DD(\Schur_{3,1}) & = \, (\Schur_3+\Schur_{2,1})\circ \Psi^1+\Schur_{1,1}\circ \Psi^2 -\Psi^4 \\
			\DD(\Schur_{2,2}) & = \, \Schur_{2,1}\circ \Psi^1 + (\Schur_2-\Schur_{1,1})\circ \Psi^2 - \Schur_1\circ \Psi^3 \\
			\DD(\Schur_{2,1,1}) & = \, (\Schur_{2,1}+\Schur_{1,1,1})\circ \Psi^1 - \Schur_2\circ \Psi^2+\Psi^4 \\
			\DD(\Schur_{1,1,1,1}) & = \, \Schur_{1,1,1}\circ \Psi^1-\Schur_{1,1}\circ \Psi^2+\Schur_1\circ \Psi^3 - \Psi^4.
		\end{align*}
		Corollary \ref{cor:plethysm formula} then yields 
		\[ \Schur_2[\Schur_2]=\Schur_4+\Schur_{2,2}\, . \]
	\end{example}
	
	\begin{example}
		Let us compute the plethysm $\Schur_2[\Schur_{1,1}]$. By a direct computation we see that
		\[ \DD(\Schur_2[\Schur_{1,1}])=(\Schur_{1,1,1}+\Schur_{2,1})\circ \Psi^1+(\Schur_2-2\Schur_{1,1})\circ \Psi^2-\Psi^4, \]
		so that we eventually get
		\[ \Schur_2[\Schur_{1,1}]=\Schur_{2,2}+\Schur_{1,1,1,1}. \]
	\end{example}
	
	\begin{example}\label{example:e2 di e2}
		Let us compute the plethysm $\Schur_{1,1}[\Schur_{1,1}]$. By a direct computation we find
		\[ \DD(\Schur_{1,1}[\Schur_{1,1}])=(\Schur_{2,1}+\Schur_{1,1,1})\circ \Psi^1-\Schur_2\circ \Psi^2+\Psi^4, \]
		from which it follows that $\Schur_{1,1}[\Schur_{1,1}]=\Schur_{2,1,1}$ is irreducible.
	\end{example}
	
	%%%%%%%%%%%%%%%%%%%%%%%%%%%%%%%%%%%%
	%%%%%%%%%%%%%%%%%%%%%%%%%%%%%%%%%%%%
	\section{An application}\label{section:applications}
	
	We start with a working definition. First of all, recall that we view Young diagrams as subsets of the south-east (henceforth, {\em positive}) quadrant in the cartesian plane.
	
	\begin{definition}
		Let $A$ be a region in the positive quadrant. We say that a symmetric function $f\in\Rep$ is \emph{supported} on $A$ if when expressed in the Schur function basis $f=\sum_\lambda m_\lambda\Schur_\lambda$, then each $\lambda$ such that $m_\lambda \neq 0$ is contained in $A$.
	\end{definition}
	
	All partitions occurring in the linear combination expressing $\D(\Schur_\lambda)$ in the Schur function basis are contained in $\lambda$. 
	This hints at the fact that if $\D(f)$ is supported in a certain region $A$, then the same should, at least approximately, hold for $f$ itself.
	However, this fails dramatically since cancellations may occur. For instance, $\D(\power_n) = n = n \Schur_{\emptyset}$ has empty support, whereas the Schur functions necessary to express $\power_n$ as a linear combination do occupy a non-trivial, indeed quite large, hook-like region (see Example \ref{example:p n come schur alpha}).
	
	The purpose of this section is to show that the above counterexamples are the only possible ones: we prove this fact for (possibly infinite) regions of {\em rectangular} form, but strongly believe some generalized statement to hold. We will then use this fact to gather inductive information on the shape of 
	Schur diagrams occurring in a plethysm.

	\begin{definition}\label{support}
		We say that $f = \sum_\lambda m_\lambda\Schur_\lambda\in\Rep$ is \emph{supported on $t$ columns} (resp.\ \emph{on $t$ rows}) if each Schur function $\Schur_\lambda$  with $m_\lambda \neq 0$ satisfies $\lambda_1\leq t$ (resp.\ $\ell(\lambda)\leq t$).
	\end{definition}

	\begin{remark}
		If $f$ is supported on $t$ columns, then it is supported on $N$ columns for every $N\geq t$. If we want to specify that $t$ is the smallest integer such that $f$ is supported on $t$ columns, then we say that $f$ is supported on \emph{exactly} $t$ columns (i.e.\ there exists at least one Schur functions $\Schur_\lambda$ in the expression of $f$ with $\lambda_1=t$). The same applies to rows.
	\end{remark}
	
	\begin{remark}\label{rmk:N minore di d}
		If $f\in\Rep^d$ is homogeneous of degree $d$, then $f$ can be uniquely written as a linear combination of Schur functions $\Schur_\lambda$ with $|\lambda|=d$. In particular, $\lambda_1\leq d$, i.e.\ $f$ is always supported on $d$ columns. Similarly, $\ell(\lambda)\leq d$, so that $f$ is always supported on $d$ rows.
		
		Notice also that if $|\lambda|=d$, then $\lambda_1=d$ if and only if $\lambda=(d)$ and $\ell(\lambda)=d$ if and only if $\lambda=(1^d)$.
	\end{remark}

	By Proposition \ref{prop:geometric interpretation of D}, if $f\in\Rep^d$ is supported on $t$ columns (resp.\ $t$ rows), then also $\D(f)$ is supported on $t$ columns (resp.\ $t$ rows). The converse is not true though, for example we know that $\D(\power_d)=d$ is supported on $0$ columns, but $\power_d$ is supported on $d$ columns.
	The next result claims that this is essentially the only exception.
	
	\begin{thm}\label{thm:supported on columns}
		Let $f\in\Rep^d$ be a homogeneous symmetric function.
		\begin{enumerate}
			\item If $\D(f)$ is supported on $t\leq d-1$ columns, then
			\[ f=g + m\power_d, \]
			with $g$ supported on $t$ columns.
			\item If $\D(f)$ is supported on $t\leq d-1$ rows, then
			\[ f=g + m\power_d, \]
			with $g$ supported on $t$ rows.
		\end{enumerate}
	\end{thm}
	
	Notice that any $f\in\Rep^d$ is supported on $d$ columns and $d$ rows (see Remark \ref{rmk:N minore di d}), so $\D(f)$ is supported on $d-1$ columns and $d-1$ rows; so the hypotheses are sharp.
	We are going to prove the theorem in Section \ref{section:proof of columns}; now we see a consequence about plethysm coefficients.
	
	\begin{corollary}\label{cor:alcuni coefficienti sono zero}
		If $\Schur_\nu$ occurs in the plethysm $\Schur_\lambda[\Schur_\mu]$, then 
		\[ \nu_1\leq |\lambda|\, \mu_1\qquad\mbox{and}\qquad\ell(\nu)\leq|\lambda|\, \ell(\mu)\,. \]
		\begin{proof}
			First of all, when $f \in \Rep^d$ is homogeneous, then all summands in $\D(f) = \sum_{n \geq 1} \power_n^\perp(f)$ have different degrees, so that we may replace $\D$ by $\DD$ in Definition \ref{support}.
			
			Let us prove the first statement, namely that $\nu_1\leq |\lambda|\, \mu_1$.
			We are going to bound the support of $\DD(\Schur_\lambda[\Schur_\mu])$ using Theorem \ref{thm:plethysm product}. By formula \eqref{eqn:D s} we have
			\[
			\DD(\Schur_\lambda)[\Schur_\mu]=\sum_{(i,j)\in\lambda}(-1)^{\varepsilon(i,j)}\,\Schur_{\lambda_{(i,j)}}[\Schur_\mu]\circ\Psi^{\hook_{i,j}(\lambda)}.
			\]
			By induction on the size of $\lambda$, we have that if $\Schur_\delta$ occurs in $\Schur_{\lambda_{(i,j)}}[\Schur_\mu]$, then $\delta_1\leq |\lambda_{(i,j)}|\, \mu_1$. Notice that $|\lambda_{(i,j)}|=|\lambda|-\hook_{i,j}(\lambda)$.
			Again by formula \eqref{eqn:D s} we have
			\[
			\DD(\Schur_\mu)=\sum_{(\ell,m)\in\mu}(-1)^{\varepsilon(\ell,m)}\,\Schur_{\mu_{(\ell,m)}}\circ\Psi^{\hook_{\ell,m}(\mu)}.
			\]
			By Theorem \ref{thm:plethysm product}, we now need to compute the product
			\[ (\Schur_\delta\circ\Psi^{\hook_{i,j}(\lambda)})\circ(\Schur_{\mu_{(\ell,m)}}\circ\Psi^{\hook_{\ell,m}(\mu)})=\Schur_{\delta}(\power_{\hook_{i,j}(\lambda)}[\Schur_{\mu_{(\ell,m)}}])\circ\Psi^{\hook_{i,j}(\lambda)\hook_{\ell,m}(\mu)}. \]
			Since the number of columns of $\mu_{(\ell,m)}$ is less than or equal to $\mu_1$, by another induction on the first row of $\mu$ we have that each factor in $\power_{\hook_{i,j}(\lambda)}[\Schur_{\mu_{(\ell,m)}}]$ is of the form $\Schur_\tau$, with $\tau_1\leq \mu_1\,\hook_{i,j}(\lambda)$.
			Finally, by Lemma \ref{elementary lemma}, if $\Schur_\sigma$ is a factor of $\Schur_\delta\,\Schur_\tau$, then
			\[
			\sigma_1\leq \delta_1+\tau_1\leq |\lambda|\, \mu_1.
			\]
			In other words, $\DD(\Schur_\lambda[\Schur_\mu])$ is supported on $|\lambda|\, \mu_1$ columns.
			By Theorem \ref{thm:supported on columns}, it then follows that 
			\begin{equation}\label{eqn:prima entrata del pletismo}
				\Schur_\lambda[\Schur_\mu]=\sum_\nu m_\nu\Schur_\nu + m\power_d,
			\end{equation}
			where each $\Schur_\nu$ is supported on $|\lambda|\, \mu_1$ columns. Here $d=|\lambda||\mu|$.
			
			Now, if $|\lambda|=1$, then there is nothing to prove, so we can assume that $|\lambda|\geq2$. Moreover, the claim is also trivial if $\mu$ has length $1$ (i.e.\ $\mu=(\mu_1)$), so we can further suppose that $|\mu|>\mu_1$. 
			Under these assumptions $|\lambda|\, \mu_1< d-1$ so that in the expression \eqref{eqn:prima entrata del pletismo} we always have the two terms $m\Schur_{(d)}-m\Schur_{(d-1,1)}$. If $m\neq0$, then there exists at least one negative term, which is in contradiction with Schur positivity of plethysms (\cite[Appendix~I.A]{Macdonald}). Therefore $m=0$ and the first statement is proved.
			
			The second statement, namely $\ell(\nu)\leq|\lambda|\, \ell(\mu)$, follows in the same way. Alternatively, one can deduce it from the first statement by passing to the conjugate partitions. In fact, if $\Schur_\lambda[\Schur_\mu]=\sum_{\nu} c^\nu_{\lambda,\mu}\,\Schur_\nu$, then by \cite[I.8, Example 1]{Macdonald} (see also \cite[Lemma~2.1]{LO23})
			\[ \Schur_{\lambda^*}[\Schur_{\mu^*}]=\sum c^{\nu^*}_{\lambda^{(*)},\mu^*}\Schur_{\nu^*}, \]
			where $\lambda^{(*)}=\lambda$ if $|\mu|$ is even and $\lambda^{(*)}=\lambda^*$ if $|\mu|$ is odd.
		\end{proof}
	\end{corollary}
	
	The following elementary result has been used in the proof.
	
	\begin{lemma}\label{elementary lemma}
		If $\Schur_\sigma$ is a factor of $\Schur_\delta\,\Schur_\tau$, then
		\[
		\sigma_1\leq \delta_1+\tau_1.
		\]
		\begin{proof}
			For any partition $\lambda$, the greatest power of the variable $x_1$ in $\Schur_\lambda(x_1,x_2,\dots)$ is $x_1^{\lambda_1}$, from which the claim follows.
		\end{proof}
	\end{lemma}
	
	\begin{remark}
		As already remarked in the Introduction, the claim of Corollary \ref{cor:alcuni coefficienti sono zero} is well-known and can be deduced, for example, from the fact that $\Schur_{\lambda}[\Schur_\mu]$ is a factor of $\Schur_\mu^{\otimes|\lambda|}$. It is worth noticing that our proof does not make use of the Littlewood--Richardson formula.
	\end{remark}
	
	%%%%%%%%%%%%%%%%%%%%%%%%%%%%%%%%%%%%%%%%%%
	\subsection{Proof of Theorem \ref{thm:supported on columns}}\label{section:proof of columns}

	Before starting the proof, let us work out by hand the case when $d=4$ (the cases $d=2$ and $d=3$ are easier and can be worked out by the reader in a similar way).
	
	\begin{example}
		Suppose that $f\in\Rep^4$, so that
		\[ f=a\,\Schur_4+b\,\Schur_{3,1}+c\,\Schur_{2,2}+d\,\Schur_{2,1,1}+e\,\Schur_{1,1,1,1}. \]
		Then 
		\begin{align*}
			\D(f) = & \, a\,(\Schur_3+\Schur_{2}+\Schur_1+1) + \\
			& \, b\,([\Schur_3+\Schur_{2,1}]+\Schur_{1,1}-1) + \\
			& \, c\,(\Schur_{2,1}+[\Schur_2-\Schur_{1,1}]-\Schur_1) + \\
			& \, d\,([\Schur_{2,1}+\Schur_{1,1,1}]-\Schur_2+1)+\\
			& \, e\,(\Schur_{1,1,1}-\Schur_{1,1}+\Schur_1-1).
		\end{align*}
		Recall that $\power_4=\Schur_4-\Schur_{3,1}+\Schur_{2,1,1}-\Schur_{1,1,1,1}$.
		\begin{itemize}
			\item If $\D(f)$ is supported on exactly $3$ columns, then there are no relations among the coefficients and we can write $f$ as
			\[ f=a\,\power_4+(b+a)\,\Schur_{3,1}+c\,\Schur_{2,2}+(d-a)\,\Schur_{2,1,1}+(e+a)\,\Schur_{1,1,1,1}. \]
			\item If $\D(f)$ is supported on exactly $2$ columns, then we get the relation $a+b=0$, so that
			\[ f=a\,\power_4 + c\,\Schur_{2,2}+(d-a)\,\Schur_{2,1,1}+(e+a)\Schur_{1,1,1,1} \]
			\item If $\D(f)$ is supported on exactly $1$ column, then in addition to $a+b=0$ we also get the relations $a+c-d=0$ and $b+c+d=0$, which give $c=0$ and $d-a=0$, i.e.
			\[ f=a\,\power_4 + (e+a)\Schur_{1,1,1,1}. \]
			\item Finally, if $\D(f)$ is supported on $0$ columns, then we get the further relation $e+a=0$ from which it follows that $f=a\,\power_4$.
		\end{itemize}
	\end{example}
	
	\proof[Proof of Theorem \ref{thm:supported on columns}]
	Let us first prove Theorem \ref{thm:supported on columns}.(1).
	Since $f\in\Rep^d$, we have that $\D(f)$ is supported on $d-1$ columns. 
	In this case, we can write 
	\[ 
	f=\sum_\lambda m_\lambda\Schur_\lambda+m\Schur_d, 
	\]
	where the sum runs over all $\Schur_\lambda$ supported on $d-1$ columns. Since $\power_d=\Schur_d-\Schur_{(d-1,1)}+\cdots+(-1)^{d-1}\Schur_{1^d}$, this can be uniquely written as
	\[ f=\sum_\lambda m'_\lambda\Schur_\lambda+m\power_d. \]
	Therefore the claim is always true if $\D(f)$ is supported on $d-1$ columns.
	
	We now proceed by reverse induction on the number of columns, using the previous case as base of the induction. Let us then suppose that the theorem has been proved for all $\tilde{f}\in\Rep^d$ such that $\D(\tilde{f})$ is supported on $d-k+1$ columns. 
	
	For any $k\geq2$, let $\D(f)$ be supported on $d-k$ columns. Then it is also supported on $d-k+1$ columns. By reverse induction we have that 
	\[ f=g+m\power_d, \]
	where $g$ is supported on $d-k+1$ columns. Then it can be written as
	\[ g=a_1\,\Schur_{d-k+1,k-1}+a_2\,\Schur_{d-k+1,k-2,1}+\cdots +a_r\,\Schur_{d-k+1,1^{k-1}}+g', \]
	with $g'$ supported on $d-k$ columns. Here $r$ is less than or equal to the number of partitions of $k-1$ in parts that are at most $d-k+1$.
	
	We have $\D(f)=md+\D(g)$, so that $\D(f)$ is supported on $d-k$ columns if and only if $\D(g)$ is supported on $d-k$ columns.
	Notice also that $\D(g')$ is supported on $d-k$ columns by construction.
	Therefore we can focus on 
	\[ h=a_1\,\Schur_{d-k+1,k-1}+a_2\,\Schur_{d-k+1,k-2,1}+\cdots +a_r\,\Schur_{d-k+1,1^{k-1}}. \]
	
	We claim that $h=0$, which will conclude the proof of the first part of the theorem.
	
	Let us denote by $\D(h)_{d-k+1}$ the part of $\D(h)$ supported on $d-k+1$ columns. By Proposition \ref{prop:geometric interpretation of D}, this is composed of those Schur functions obtained from the ones in $h$ by removing hooks centered in boxes away from the first row. Let us make this observation more formal by introducing some notations. If $\lambda=(\lambda_1,\dots,\lambda_n)$ is a partition, we put $\underline{\lambda}=(\lambda_2,\dots,\lambda_n)$ for the partition obtained by removing the first row from $\lambda$. Similarly, if $N\geq\lambda_1$, we put $\hat\lambda^N=(N,\lambda_1\dots,\lambda_n)$ for the partition obtained by adding a row to $\lambda$. For example, $\hat{\underline\lambda}^{\lambda_1}=\lambda$.
	With this notation at hand, if $(i,j)\in\underline\lambda$ is a box, then 
	\[ \lambda_{(i+1,j)}=\hat{\underline\lambda}_{(i,j)}^{\lambda_1},
	\]
	i.e.\ we first remove the first row from $\lambda$, remove the hook centred in $(i,j)$ from $\underline\lambda$, and finally we add back the first row.  Finally, if $F$ is any expression in Schur functions $\Schur_\lambda$, then we denote by $\underline{F}$ and $\hat{F}^N$, respectively, the corresponding expressions in $\Schur_{\underline\lambda}$ and $\Schur_{\hat\lambda^N}$.
	
	Let us go back to the proof. By the discussion above and Proposition \ref{prop:geometric interpretation of D}, we have 
	\[ \D(h)_{d-k+1}=\widehat{\D(\underline{h})}^{d-k+1}. \]
	
	By hypothesis we have that $ \D(h)_{d-k+1}=0$, which translates to $\D(\underline{h})=0$. On the other hand, by Theorem \ref{thm:quasi isometria} we know that $\D$ is injective on homogeneous components of positive degree. Since $k\geq2$ by assumption, we must have $\underline{h}=0$, which implies $h=0$, thus proving the claim.
	
	Finally, Theorem \ref{thm:supported on columns}.(2) follows from the first part by passing to the conjugate partitions. In fact, there exists an involution $w\colon\Rep\to\Rep$ such that $w(\Schur_\lambda)=\Schur_{\lambda^*}$ (see e.g.\ \cite[I.(2.7)]{Macdonald}) and clearly $f\in\Rep^d$ is supported on $t$ rows if and only if $w(f)$ is supported on $t$ columns (see also Example \ref{example:D of transpose} for the composition $\D\circ w$).
	\endproof

	%%%%%%%%%%%%%%%%%%%%%%%%%%%%%%%%%%%
	%%%%%%%%%%%%%%%%%%%%%%%%%%%%%%%%%%%
	%%%%%%%%%%%%%%%%%%%%%%%%%%%%%%%%%%%


\begin{thebibliography}{COSSZ24}
		
		\bibitem[Bof91]{Boffi}
		G.~Boffi.
		\emph{On some plethysms}.
		Adv.\ Math.\ \textbf{89}, no.~2 (1991), 107--126.
		\url{https://doi.org/10.1016/0001-8708(91)90075-I}.
		
		\bibitem[COSSZ24]{MysteryOfPlethysms}
		L.~Colmenarejo, R.~Orellana, F.~Saliola, A.~Schilling, and M.~Zabrocki.
		\emph{The mystery of plethysm coefficients}.
		Proc.\ Sympos.\ Pure Math.\ \textbf{110} (2024), 275--292.
		\url{https://doi.org/10.1090/pspum/110/02018}.
		
		\bibitem[DFO25a]{Casimiri}
		A.~D'Andrea, E.~Fatighenti, and C.~Onorati.
		\emph{Higher discriminants of vector bundles and Schur functors}.
		To appear in Ann.\ Sc.\ Norm.\ Super.\ Pisa Cl.\ Sci.;
		arXiv:2503.15365 (2025).
		\url{https://doi.org/10.2422/2036-2145.202505_020}.
		
		\bibitem[DFO25b]{DFO:Chern}
		A.~D'Andrea, E.~Fatighenti, and C.~Onorati.
		\emph{Chern character of Schur bundles}.
		arXiv:2512.11657 (2025).
		
		\bibitem[dBPW21]{dBPW21}
		M.~De Boeck, R.~Paget, and M.~Wildon.
		\emph{Plethysms of symmetric functions and highest weight representations}.
		Trans.\ Amer.\ Math.\ Soc.\ \textbf{374}, no.~11 (2021), 8013--8043.
		\url{https://doi.org/10.1090/tran/8481}.
		
		
		\bibitem[EPW14]{EPW}
		A.~Evseev, R.~Paget, and M.~Wildon.
		\emph{Character deflations and a generalization of the
			Murnaghan--Nakayama rule}.
		J.\ Group Theory \textbf{17}, no.~6 (2014), 1035--1070.
		
		\bibitem[FI20]{FI}
		N.~Fischer and C.~Ikenmeyer.
		\emph{The computational complexity of plethysm coefficients}.
		Comput.\ Complex.\ \textbf{29} (2020), article~8.
		\url{https://doi.org/10.1007/s00037-020-00198-4}.
		
		\bibitem[KM16]{KM16}
		T.~Kahle and M.~Micha{\l}ek.
		\emph{Plethysm and lattice point counting}.
		Found.\ Comput.\ Math.\ \textbf{16} (2016), 1241--1261.
		\url{https://doi.org/10.1007/s10208-015-9275-7}.
		
		\bibitem[KM18]{KM18}
		T.~Kahle and M.~Micha{\l}ek.
		\emph{Obstructions to combinatorial formulas for plethysm}.
		Electron.\ J.\ Combin.\ \textbf{25}, no.~1 (2018), Paper~P1.41.
		\url{https://doi.org/10.37236/6597}.
		
		\bibitem[LO23]{LO23}
		S.~Law and Y.~Okitani.
		\emph{Some stable plethysms}.
		Proc.\ Amer.\ Math.\ Soc.\ \textbf{151}, no.~11 (2023), 4557--4564.
		\url{https://doi.org/10.1090/proc/16556}.
		
		\bibitem[Lit36]{Littlewood}
		D.~E.~Littlewood.
		\emph{Polynomial concomitants and invariant matrices}.
		J.\ London Math.\ Soc.\ \textbf{11} (1936), 49--55.
		\url{https://doi.org/10.1112/jlms/s1-11.1.49}.
		
		\bibitem[LR11]{LoehrRemmel2011}
		N.~A.~Loehr and J.~B.~Remmel.
		\emph{A computational and combinatorial expos\'e of plethystic calculus}.
		J.\ Algebraic Combin.\ \textbf{33} (2011), 163--198.
		\url{https://doi.org/10.1007/s10801-010-0238-4}.
		
		\bibitem[Mac95]{Macdonald}
		I.~G.~Macdonald.
		\emph{Symmetric functions and Hall polynomials}.
		2nd ed., Oxford Mathematical Monographs.
		Oxford University Press, Oxford, 1995.
		
		\bibitem[PW16]{PW16}
		R.~Paget and M.~Wildon.
		\emph{Minimal and maximal constituents of twisted Foulkes characters}.
		J.\ London Math.\ Soc.\ \textbf{93}, no.~2 (2016), 301--318.
		\url{https://doi.org/10.1112/jlms/jdv070}.
		
		\bibitem[Sta]{StanleyII}
		R.~P.~Stanley.
		\emph{Enumerative combinatorics. Vol.~2}.
		Cambridge Studies in Advanced Mathematics, vol.~62.
		Cambridge University Press, Cambridge, 1999.
		With a foreword by Gian-Carlo Rota and Appendix~1 by Sergey Fomin.
		
		\bibitem[Sta00]{Stanley}
		R.~P.~Stanley.
		\emph{Positivity problems and conjectures in algebraic combinatorics}.
		In \emph{Mathematics: Frontiers and Perspectives},
		American Mathematical Society, Providence, RI, 2000, 295--319.
		
		\bibitem[Yan16]{Yanagida}
		S.~Yanagida.
		\emph{Whittaker vector of deformed Virasoro algebra and Macdonald
			symmetric functions}.
		Lett.\ Math.\ Phys.\ \textbf{106} (2016), 395--431.
		\url{https://doi.org/10.1007/s11005-016-0821-2}.
		
		\bibitem[Zab13]{Zabrocki}
		M.~Zabrocki.
		\emph{Introduction to symmetric functions}.
		Lecture notes, 2013, 111~pp.
		\url{https://garsia.math.yorku.ca/ghana03/mainfile.pdf}.
		
	\end{thebibliography}
\end{document}